\documentclass{amsart}

\usepackage{fancyhdr}
\usepackage{lastpage}
\usepackage{stmaryrd,yhmath}

\newcommand{\bdis}{\begin{displaymath}}
\newcommand{\edis}{\end{displaymath}}
\newcommand{\be}{\begin{equation}}
\newcommand{\ee}{\end{equation}}
\newcommand{\mbb}{\mathbb}
\newcommand{\mcal}{\mathcal}

\newcommand{\vp}{\varphi}

\theoremstyle{definition}

\newtheorem{cor}[]{Corollary}

\theoremstyle{remark}
\newtheorem{remark}[]{Remark}

\newtheorem*{mydef1}{{\bf Theorem}}

\newtheorem*{mydef5}{{\bf Lemma}}

\numberwithin{equation}{section}

\begin{document}

\title{Jacob's ladders, the iterations of Jacob's ladder $\varphi^k_1(t)$ and asymptotic formulae for the integrals of the products
$Z^2[\varphi^n_1(t)]Z^2[\varphi^{n-1}(t)]\cdots Z^2[\varphi^0_1(t)]$ for arbitrary fixed $n\in\mbb{N}$}

\author{Jan Moser}

\address{Department of Mathematical Analysis and Numerical Mathematics, Comenius University, Mlynska Dolina M105, 842 48 Bratislava, SLOVAKIA}

\email{jan.mozer@fmph.uniba.sk}

\keywords{Riemann zeta-function}

\begin{abstract}
In this paper we introduce the iterations $\varphi^k_1(t)$ of the Jacob's ladder. It is proved, for example, that the mean-value of the
product $$Z^2[\varphi^n_1(t)]Z^2[\varphi^{n-1}(t)]\cdots Z^2[\varphi^0_1(t)]$$ over the segment $[T,T+U]$ is asymptotically equal to
$\ln^{n+1}T$. Nor the case $n=1$ cannot be obtained in known theories of Balasubramanian, Heath-Brown and Ivic.
\end{abstract}

\maketitle

\section{Results}

Let
\begin{eqnarray} \label{1.1e}
& &
y=\frac{1}{2}\vp(t)=\vp_1(t);\ \vp^0_1(t)=t;\ \vp_1^1(t)=\vp_1(t),\ \vp_1^2(t)=\vp_1(\vp_1(t)),\ \nonumber \\
& & \ \\
& &
\dots ,\ \vp^k_1(t)=\vp_1(\vp_1(\dots (\vp_1(t))\dots ),\ t\in [T,T+U] , \nonumber
\end{eqnarray}
where $\vp_1^k(t)$ denotes the $k$th iteration of the Jacob's ladder $y=\vp(t),\ t\geq T_0[\vp_1]$. The following Theorem holds true.

\begin{mydef1}
Let
\be \label{1.2e}
T\geq T_{00}[\vp_1,n]=\max\{ 2T_0[\vp_1],e^{2(n+1)}\},\ U=T^{1/3+2\epsilon} .
\ee
Then for every fixed $n\in\mbb{N}$ the following is true
\be \label{1.3e}
\int_T^{T+U}\prod_{k=0}^n Z^2[\vp^k_1(t)]{\rm d}t\sim U\ln^{n+1}T,\ T\to\infty ,
\ee
where
\begin{equation}
\vp_1^k(t)\geq T_0[\vp_1[,\ k=0,1,\dots ,n+1,\ t\in [T,T+U], \tag{A}
\end{equation}
\begin{equation}
\vp_1^k(T+U)-\vp_1^k(t)\sim U,\ k=0,1,\dots ,n+1, \tag{B}
\end{equation}
\begin{equation}
\vp_1^{k-1}(T)-\vp_1^k(T+U)\sim (1-c)\frac{T}{\ln T},\ k=0,1,\dots ,n, \tag{C}
\end{equation}
\begin{equation}
\rho\left\{ \left[ \vp_1^{k-1}(T),\vp_1^{k-1}(T+U)\right];\left[ \vp_1^k(T),\vp_1^k(T+U)\right]\right\}\sim (1-c)\frac{T}{\ln T}, \tag{D}
\end{equation}
and $\rho$ denotes the distance of the corresponding segments.
\end{mydef1}

\begin{remark}
The system of the iterated segments
\bdis
\left[ \vp_1^n(T),\vp_1^n(T+U)\right],\ \left[ \vp_1^{n-1}(T),\vp_1^{n-1}(T+U)\right], \dots , [T,T+U]
\edis
is the disconnected set of segments distributed from right to left (see (C)) and the neighbouring segments unboundedly recede each from other
(see (D), $\rho\to \infty$, as $T\to \infty$, comp. \cite{6}, Remark 3).
\end{remark}

\begin{remark}
Let us mention the formula (1.3) especially for the prime numbers of Fermat-Gauss $n=17,257,65 537$ and for the Skewes constant
\bdis
n=10^{10^{10^{34}}} .
\edis
It is obvious that nor the formula ($n=2$)
\bdis
\int_T^{T+U}Z^2[\vp_1(\vp_1(t))]Z^2[\vp_1(t)]Z^2(t){\rm d}t\sim U\ln^3 T
\edis
cannot be reached in known theories of Balasubramanian, Heath-Brown and Ivic (see \cite{1}).
\end{remark}

This paper is a continuation of the series of papers \cite{2}-\cite{7}.

\section{Consequences of the Theorem}

Using the mean-value theorem in (1.3) we obtain
\begin{cor}
\begin{eqnarray} \label{2.1e}
& &
\prod_{k=0}^nZ^2[\vp_1^k(\tau)]\sim\ln^{n+1}T,\\
& &
T\to\infty, \vp_1^k(\tau)\in(\vp_1^k(T),\vp_1^k(T+U)), k=0,1,\dots ,n, \tau=\tau(T,n). \nonumber
\end{eqnarray}
\end{cor}
From (2.1) we obtain
\begin{cor}
\be \label{2.2e}
\prod_{k=0}^n \left| Z[\vp_1^k(\tau)]\right|^{\frac{2}{n+1}}\sim \ln T,
\ee
\be \label{2.3e}
\frac{1}{n+1}\sum_{k=0}^n\ln|Z[\vp_1^k(\tau)]|\sim \frac{1}{2}\ln\ln T .
\ee
\end{cor}

Next, by the known inequalities for harmonic, geometric and arithmetic means we have
\begin{cor}
\be \label{2.4e}
(1-\epsilon)\ln T\leq \frac{1}{n+1}\sum_{k=0}^n |Z[\vp_1^k(\tau)]|^{\frac{2}{n+1}} ,
\ee
\be\label{2.5e}
\frac{1}{(1+\epsilon)\ln T}< \frac{1}{n+1}\sum_{k=0}^n|Z[\vp_1^k(\tau)]|^{-\frac{2}{n+1}} .
\ee
\end{cor}

\begin{remark}
Some new type of the nonlocal interaction of the values
\bdis
\left\{ Z^2[\vp_1^k(t)]\right\}_{k=0}^n
\edis
of the signal
\bdis
Z(t)=e^{i\vartheta(t)}\zeta\left(\frac{1}{2}+it\right)
\edis
over the system of disconnected segments
\bdis
\bigcup_{k=0}^n \left[ \vp_1^k(T),\vp_1^k(T+U)\right]
\edis
is expressed by formulae (1.3), (2.1)-(2.5) for the iterated Jacob's ladder $\vp_1^k(t)$.
\end{remark}

\section{Lemma}

We start with the formula (see \cite{2}, (3.5), (3.9))
\bdis
Z^2(t)=\Phi^\prime_\vp\left[\vp(t)\right]\frac{{\rm d}\vp(t)}{{\rm d}t},\ t\geq T_0[\vp] ,
\edis
where (see \cite{4}, (1.5))
\bdis
\Phi^\prime_\vp[\vp(t)]=\frac{1}{2}\left\{ 1+\mcal{O}\left(\frac{\ln\ln t}{\ln t}\right)\right\}\ln t .
\edis
Next we have
\be \label{3.1e}
\tilde{Z}^2(t)=\frac{{\rm d}\vp_1(t)}{{\rm d}t},\ t\in [T,T+U],\ U\in\left(\left.  0,\frac{T}{\ln T}\right]\right. ,
\ee
by (1.1) we have
\be \label{3.2e}
\tilde{Z}^2(t)=\frac{Z^2(t)}{2\Phi^\prime_\vp[\vp(t)]}=\frac{Z^2(t)}{\left\{ 1+\mcal{O}\left(\frac{\ln\ln t}{\ln t}\right)\right\}\ln t} .
\ee
Then we obtain from (3.1) the following lemma (comp. \cite{6}, (2.5)).

\begin{mydef5}
For every integrable function (in the Lebesgue sense) $f(x), x\in [\vp_1(T),\vp_1(T+U)]$ the following is true
\be\label{3.3e}
\int_T^{T+U} f[\vp_1(t)]\tilde{Z}^2(t){\rm d}t=\int_{\vp_1(T)}^{\vp_1(T+U)}f(x){\rm d}x,\ U\in \left(\left.  0,\frac{T}{\ln T}\right]\right.  .
\ee
\end{mydef5}

\section{Proof of the Theorem}

\subsection{}

From the formula (see (1.1), (1.2), \cite{2}, (6.2))
\be \label{4.1e}
t-\vp_1^1(t)\sim (1-c)\pi(t)\sim (1-c)\frac{t}{\ln t}
\ee
we have
\begin{eqnarray} \label{4.2e}
& &
\vp_1^1(t)-\vp_1^2(t)\sim (1-c)\frac{t}{\ln t} ,\nonumber \\
& &
\vp_1^2(t)-\vp_1^3(t)\sim (1-c)\frac{t}{\ln t}, \nonumber \\
& & \vdots \\
& &
\vp_1^n(t)-\vp_1^{n+1}(t)\sim (1-c)\frac{t}{\ln t}, \nonumber
\end{eqnarray}
and by (4.1) (4.2) we obtain
\bdis
\vp_1^{n+1}(t)>t\left\{ 1-\frac{(1+\epsilon)(1-c)(n+1)}{\ln t}\right\}\geq T\left( 1-\frac{n+1}{\ln T}\right)\geq \frac{1}{2}T\geq T_0[\vp_1],
\edis
i.e. (A).

\subsection{}

By comparison of the formula (see \cite{3}, (1.5))
\bdis
\int_T^{T+U}Z^2(t){\rm d}t\sim U\ln T,\ U=T^{1/3+2\epsilon} ,
\edis
where $1/3$ is the exponent of Balasubramanian, and our formula
\bdis
\int_T^{T+U} Z^2(t){\rm d}t\sim \left\{ \vp_1(T+U)-\vp_1(T)\right\}\ln T ,
\edis
(see \cite{4}, (1.2)) we have
\bdis
\vp_1(T+U)-\vp_1(T)\sim U ,
\edis
by comparison in the cases $T\to \vp_1^1(T),\ T+U\to \vp_1^1(T+U),\ \dots$ we obtain
\begin{eqnarray*}
\vp_1^2(T+U)-\vp_1^2(T)&\sim & \vp_1^1(T+U)-\vp_1^1(T) , \\
& \vdots & \\
\vp_1^{n+1}(T+U)-\vp_1^{n+1}(T) & \sim & \vp_1^n(T+U)-\vp_1^n(T) ,
\end{eqnarray*}
i.e. (B).

\subsection{}

By (4.2), $t\to T$ we have
\bdis
\vp_1^1(T)-\vp_1^2(T)\sim (1-c)\frac{T}{\ln T} ,
\edis
i.e.
\be \label{4.3e}
\vp_1^1(T)-\vp_1^2(T+U)+\{ \vp_1^2(T+U)-\vp_1^2(T)\}\sim (1-c)\frac{T}{\ln T} ,
\ee
since (see (B))
\bdis
\vp_1^2(T+U)-\vp_1^2(T)\sim U=T^{1/3+2\epsilon}
\edis
then from (4.3) the asymptotic formula
\bdis
\vp_1^1(T)-\vp_1^2(T+U)\sim (1-c)\frac{T}{\ln T}
\edis
follows. Similarly we obtain all asymptotic formulae in (C). The proposition (D) follows from (C).

\subsection{}

From (3.1) by (3.3) we have
\begin{eqnarray*}
& &
\int_T^{T+U}\prod_{k=0}^n \tilde{Z}^2[\vp_1^k(t)]{\rm d}t=\int_T^{T+U}\prod_{k=0}^n \tilde{Z}^2[\vp_1^k(t)]\tilde{Z}^2(t){\rm d}t=\\
& &
=\int_T^{T+U}\prod_{k=1}^n\tilde{Z}^2[\vp_1^{k-1}(\vp_1(t))]{\rm d}\vp_1(t)=\int_{\vp_1^1(T)}^{\vp_1^1(T+U)}\prod_{k=1}^n\tilde{Z}^2
[\vp_1^{k-1}(w_1)]{\rm d}w_1=\\
& &
=\int_{\vp_1^1(T)}^{\vp_1^1(T+U)}\prod_{k=2}^n\tilde{Z}^2[\vp_1^{k-2}(\vp_1(w_1))]\tilde{Z}^2(w_1){\rm d}w_1= \\
& &
=\int_{\vp_1^2(T)}^{\vp_1^2(T+U)}\prod_{k=2}^n \tilde{Z}^2[\vp_1^{k-2}(w_2)]{\rm d}w_2=\dots = \\
& &
=\int_{\vp_1^n(T)}^{\vp_1^n(T+U)}\tilde{Z}^2[w_n]{\rm d}w_n=\vp_1^{n+1}(T+U)-\vp_1^{n+1}(T) ,
\end{eqnarray*}
i.e. the following asymptotic formula (see (B))
\be \label{4.4e}
\int_T^{T+U}\prod_{k=0}^n\tilde{Z}^2[\vp_1^k(t)]{\rm d}t=\vp_1^{n+1}(T+U)-\vp_1^{n+1}(T)\sim U
\ee
holds true. Then, from (4.4) by mean-value theorem (see (3.1), (3.2), (4.1), (4.2); $\ln \vp_1^k(t)\sim \ln t$) the formula (1.3) follows.

\thanks{I would like to thank Michal Demetrian for helping me with the electronic version of this work.}


\begin{thebibliography}{19}
%
%
%
\bibitem{1}
A. Ivic, `The Riemann zeta-function', A Willey-Interscience Publication, New York, 1985.
%
\bibitem{2}
J. Moser, `Jacob's ladders and the almost exact asymptotic representation of the Hardy-Littlewood integral', (2008).
%
\bibitem{3}
J. Moser, `Jacob's ladders and the tangent law for short parts of the Hardy-Littlewood integral', (2009).
%
\bibitem{4}
J. Moser, `Jacob's ladders and the multiplicative asymptotic formula for short and microscopic parts of the Hardy-Littlewood integral', (2009).
%
\bibitem{5}
J. Moser, `Jacob's ladders and the quantization of the Hardy-Littlewood integral', (2009).
%
\bibitem{6}
J. Moser,
`Jacob's ladders and the first asymptotic formula for the expression of the sixth order $|\zeta(1/2+i\varphi(t)/2)|^4|\zeta(1/2+it)|^2$', (2009).
%
\bibitem{7}
J. Moser,
`Jacob's ladders and the first asymptotic formula for the expression of the fifth order
$Z[\varphi(t)/2+\rho_1]Z[\varphi(t)/2+\rho_2]Z[\varphi(t)/2+\rho_3]\hat{Z}^2(t)$ for the collection of disconnected sets`, (2009).
%
\bibitem{8}
E.C. Titchmarsh,
`The theory of the Riemann zeta-function`, Clarendon Press, Oxford, 1951.


\end{thebibliography}
\end{document}